# Compressive Sensing for Polyharmonic Subdivision Wavelets With Applications to Image Analysis

Ognyan Kounchev, Damyan Kalaglarsky

**Abstract:** We apply successfully the Compressive Sensing approach for Image Analysis using the new family of Polyharmonic Subdivision wavelets. We show that this approach provides a very efficient recovery of the images based on fewer samples than the traditional Shannon-Nyquist paradigm. We provide the results of experiments with PHSD wavelets and Daubechies wavelets, for the Lena image and astronomical images.

## 1 Introduction

Compressive sensing (CS) is a relatively new approach to finding fast solutions of problems which are a priori known to be sparse. CS may be considered as a new sampling techniques which differs from the standard Shannon-Nyquist sampling approach. This concept has been developed in the context of Signal and Image Processing in the seminal papers [2], [5], [3]. In the course of development it has appeared that the basic principles and techniques of CS have been discovered at least decade before in other areas of mathematics. However CS has given a lot of results and even more hopes to improve essentially the quality and the speed of the Signal and Image Analysis.

In the present research we apply the CS methodology to a new class of multivariate wavelets, the so-called Polyharmonic Subdivision (PHSD) wavelets, [6], [7], [9], [11]. Let us recall that in [11] we have seen an amazing advantage of the PHSD combined with Vector Quantization of the coefficients over other wavelet methods of compression.

## 2 Polyharmonic subdivision wavelets

Let us recall shortly the PHSD wavelets, [7], [9], [10]. They are based on the Fourier transform in one direction $y \in \mathsf{R}^{n-1}$, and on $1d$ wavelet analysis in the direction $t$. They have demonstrated certain advantage over other established wavelet methods.

The PHSD wavelets consider the image $f(t,y)$ with $y \in \mathsf{R}^{n-1}$ and $t \in \mathsf{R}$ in the Fourier domain

$$f(t,y) = \sum_\xi f_\xi(t) e^{i\xi y},$$

where $f_\xi(t)$ is the Fourier transform in direction $y$; we assume for simplicity that the image $f$ is $2\pi-$ periodic in $y$. Then for every index $\xi$ we apply a Wavelet Analysis with wavelets $\psi_j^\xi(t)$ which depend on the index $\xi$ and on the level $j \in \mathsf{Z}$. The dependence on the index $j$ shows that this Wavelet Analysis is **non-stationary** and is based on a subdivision algorithm which is characterized by preserving the exponential functions $\left\{ t^j e^{-|\xi|t}, t^j e^{|\xi|t} \right\}_{j=0}^{p-1}$ for some integer $p$. Hence, we obtain

$$f_\xi(t) = \sum_{k,j} \psi_j^\xi(t-k) g_{k,j}^\xi.$$

This theory has been developed as a direct generalization of the famous Daubechies wavelets which are obtained in the subdivision context studied by Deslaurier and Dubuc, cf. further details [9], [10], [6], [7].

## 3 The general scheme of CS for PHSD Wavelets

The general scheme of CS is applied to a signal $x \in \mathsf{R}^N$, where $N$ is large. We assume that the signal $x$ is "sparse" in some orthonormal basis $\Phi$, and $x = \Phi\alpha$; here $\alpha$ are the coefficients of the expansion, and $\alpha = \Phi^T x$. "Sparseness" means that most entries of the coefficient vector $\alpha$ are close to zero, i.e. only few of the coefficients $\alpha_j$ are significant. Hence, the signal $x$ is "compressible" since $x' = \Phi\alpha'$ is a good approximation of $x$, where $\alpha'$ has the same coefficients as $\alpha$ but those of small size are put equal to zero. We assume that only $K$ coefficients of $\alpha$ are significant, i.e. $K \ll N$.

The main gist of CS is that the recovery of the signal $x$ is through some measurements with vectors $\theta_m \in \mathsf{R}^N$,

$$y_m = \langle x, \theta_m \rangle, \qquad m = 1, 2, \ldots, M;$$

here the measurement matrix $\Theta = \{\theta_m\}_{m=1}^M$ is properly chosen and $M$ is considerably smaller than $N$, i.e. $M \ll N$. Then the approximation to $\alpha$ is obtained by solving the problem

$$\min_\beta \|\beta\|_{\ell_1} \text{ where } y = \Theta\Phi\beta.$$

Here $\|\beta\|_{\ell_1} = \sum_j |\beta_j|$. For the details we refer to [2], [3], and the recent survey [4]. The CS has been applied also in the curvelet basis [13].

In the case of the PHSD wavelets we have the following details of the implementation. We refer also to the previous research in [10], [11], where complete details of the applications of PHSD wavelets to Image Analysis are available. The image $f(t, y)$ is a matrix (function) in two discrete variables $t$ and $y$. Respectively, we have measuring matrices $H_{t,y}^m$ playing the role of the measurement vectors $\theta_m$ above. We have the measurements

$$Y_m = \sum_{t,y} H_{t,y}^m f(t, y) \qquad \text{for } m = 1, 2, \ldots, M.$$

First of all, in the Fourier transform in direction $y$ we have

$$f(t, y) = \sum_\xi f_\xi(t) e^{i\xi y},$$

and the PHSD wavelet representation gives

$$f_\xi(t) = \sum_{k,j} \phi_j^\xi(t - k) g_{k,j}^\xi.$$

Here the basis functions $\phi^\xi$ are not only mother but may be also father wavelets (scaling functions), depending on the level of the wavelet decomposition, [10], [11].

Hence, we obtain for the measurements:

$$Y_m = \sum_{t,y} H_{t,y}^m f(t, y) = \sum_{t,y} H_{t,y}^m \left( \sum_\xi f_\xi(t) e^{i\xi y} \right)$$

$$= \sum_{t,y} H_{t,y}^m \left( \sum_{\xi} \sum_{k,j} \phi^\xi(2^j t - k) g_{k,j}^\xi e^{i\xi y} \right)$$

$$= \sum_{\xi} \sum_{k,j} g_{k,j}^\xi \left( \sum_{t,y} H_{t,y}^m \phi^\xi(2^j t - k) e^{i\xi y} \right) = \sum_{\xi} \sum_{k,j} g_{k,j}^\xi G_{\xi;t,y}^m$$

for $m = 1, 2, \ldots, M$. Thus the **CS consists in** solving the problem

$$\min \|\beta\|_{\ell_1} \text{ for } Y_m = \sum_{\xi} \sum_{k,j} \beta_{k,j}^\xi G_{\xi;t,y}^m \qquad \text{for } m = 1, 2, \ldots, M. \tag{1}$$

A slightly different approach would be to work directly in the Fourier domain by considering the measurements:

$$Y_m = \sum_{t,y} H_{t,\xi}^m f_\xi(t) \qquad \text{for } m = 1, 2, \ldots, M.$$

Hence, we obtain

$$Y_m = \sum_{t,\xi} H_{t,\xi}^m f_\xi(t) = \sum_{t,\xi} H_{t,\xi}^m \left( \sum_{k,j} \phi^\xi(2^j t - k) g_{k,j}^\xi \right)$$

$$= \sum_{k,j,\xi} g_{k,j}^\xi \left( \sum_t H_{t,\xi}^m \phi^\xi(2^j t - k) \right) = \sum_{k,j,\xi} g_{k,j}^\xi G_{k,j;\xi}^m$$

for $m = 1, 2, \ldots, M$.

Here **CS consists in** solving the minimization problem

$$\min \|\beta\|_{\ell_1} \text{ for } Y_m = \sum_{k,j,\xi} \beta_{k,j}^\xi G_{k,j,\xi}^m \qquad \text{for } m = 1, 2, \ldots, M \tag{2}$$

## 4 Experiments

In our experiments we choose different seminal images $f(t, y)$. We choose the measurement matrix $\Theta$ as in the setting (2). We make direct use of the implementation proposed in the monograph [14], and the package **splitting solvers**, available at the address:
http://thames.cs.rhul.ac.uk/126fionn/Sparse_Signal_Recipes/
Software_and_Images_files/SplittingSolvers.tar.gz

### 4.1 Lena image

This is an experiment with the Lena image (256 times 256 pixels), with 10 iterations, carried out with the **PHSD wavelets** by applying the Basis Pursuit and the Lasso methods; we consider the Sensing operator (matrix) which is represented by 50 radial lines in the Fourier domain (the Fourier transform is the measurement matrix, where we choose lines and the points on them in the Fourier domain, the number of measurements is about $M = 5020$). We see that the PSNR is almost the same with both methods; the parameter of the Lasso is $\mu = 1$ and for BP is $\gamma = 100$.

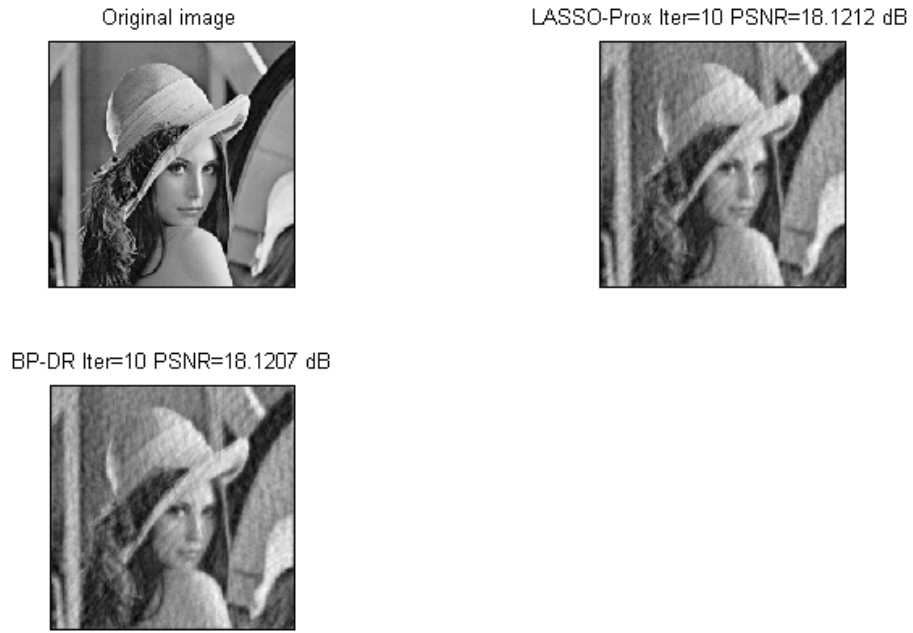

Figure 1

We compare this with the performance below of the Daubechies wavelets, where we see some advantage of the PHSD wavelets in terms of the PSNR.

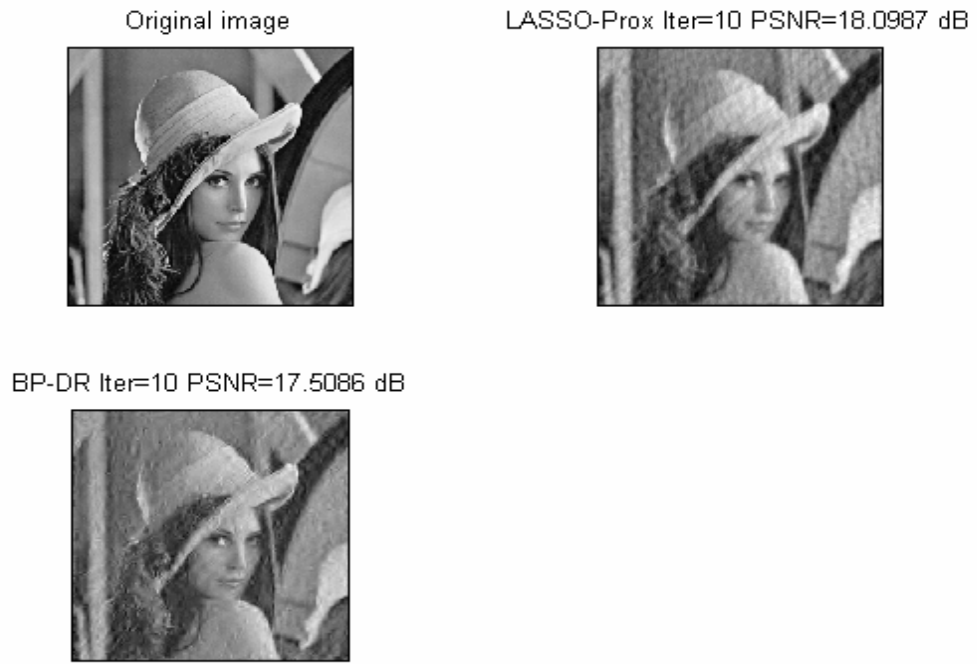

Figure 2

The next Figure contains the similar result as above but with $M = 15000$, since the measurement matrix is taken with 100 lines in the Fourier domain:

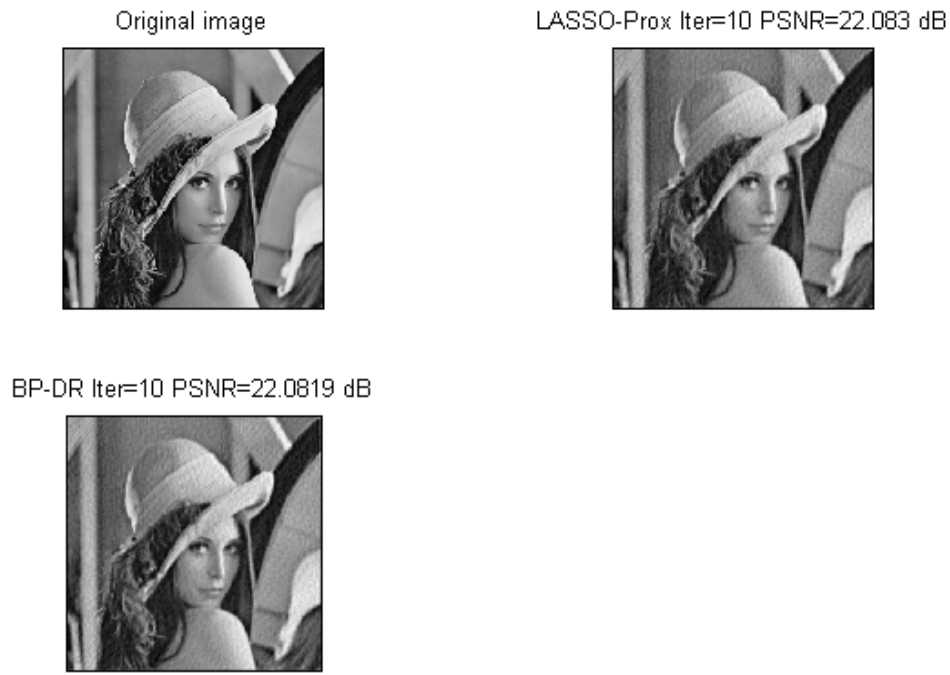

Figure 3

## 4.2 Astronomical images

We test some Astronomical images.

### 4.2.1 Chain image from Konkoly

We have chosen a **chain image** from the Konkoly plate archive. The details of the algorithms are similar to the above for the Lena image. First, we present the performance of the **PHSD wavelets**; we make 10 iterations with the Lasso and BP with the same parameters $\mu$ and $\gamma$; we take 50 lines in the Fourier domain.

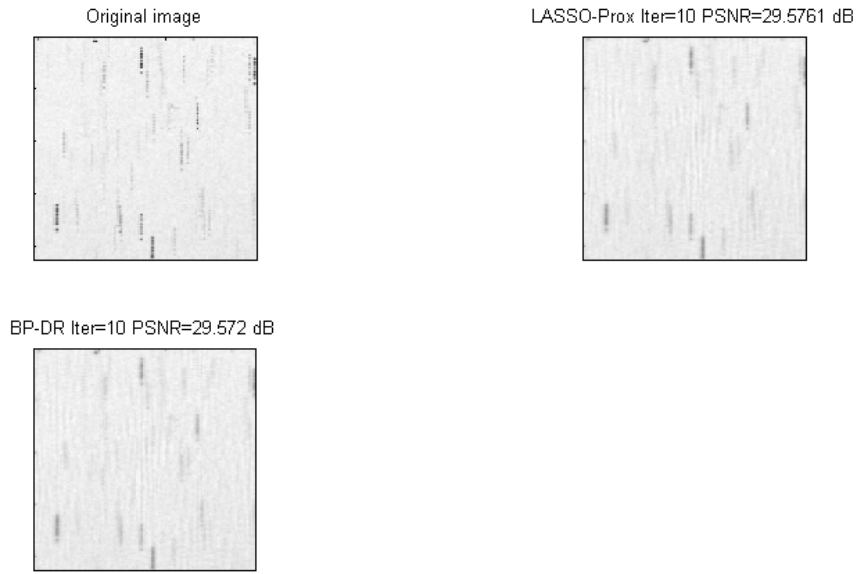

Figure 4

Next we compare it with the result of the **Daubechies wavelets** again with 50 lines of measurements in the Fourier domain:

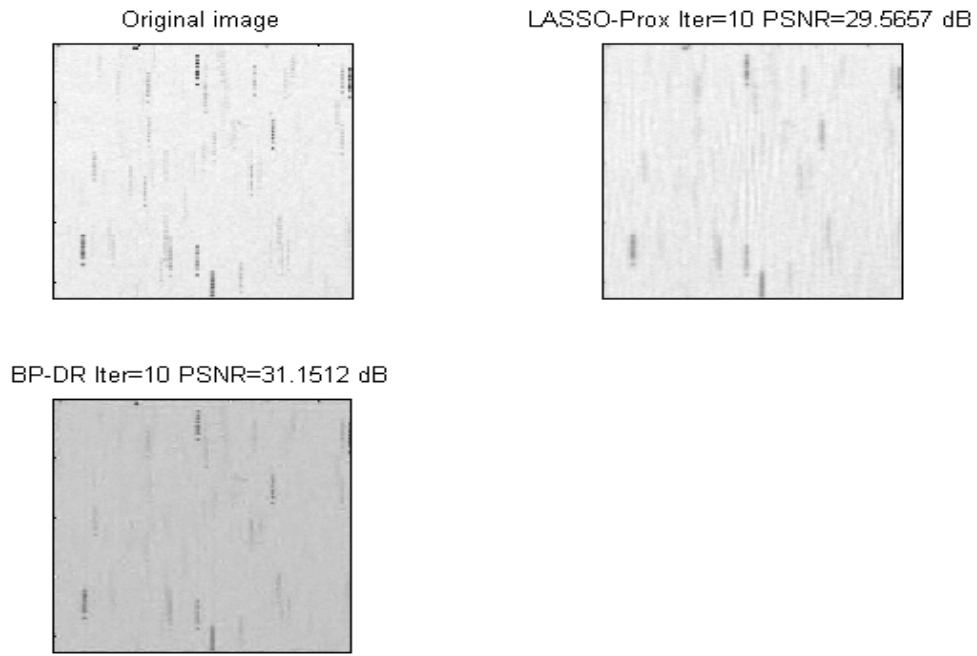

Figure 5

We see that the quality measured by PSNR is more or less the same. Similarly, we take a bigger measurement matrix (with 100 lines in the Fourier domain) and again we take **PHSD wavelets**:

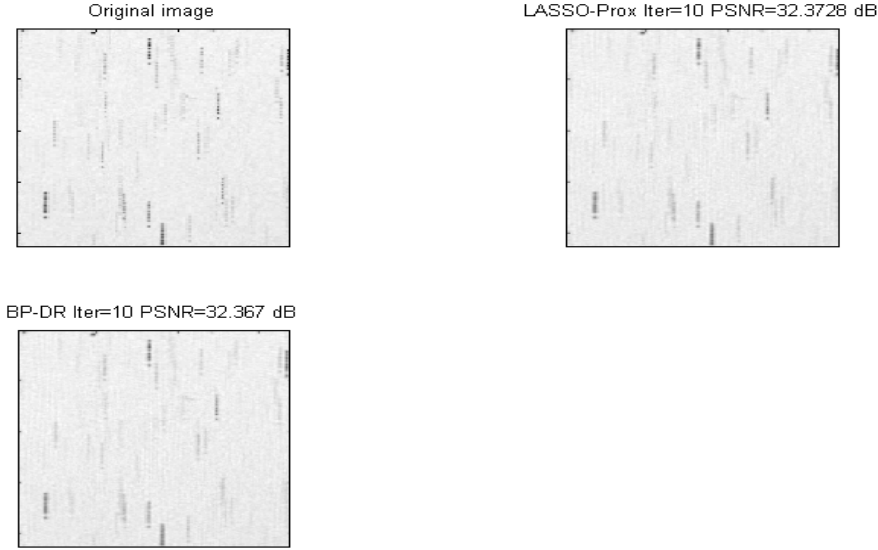

Figure 6

We compare it with the **Daubechies wavelets** with the same measurement matrices (100 lines in the Fourier domain) and this gives:

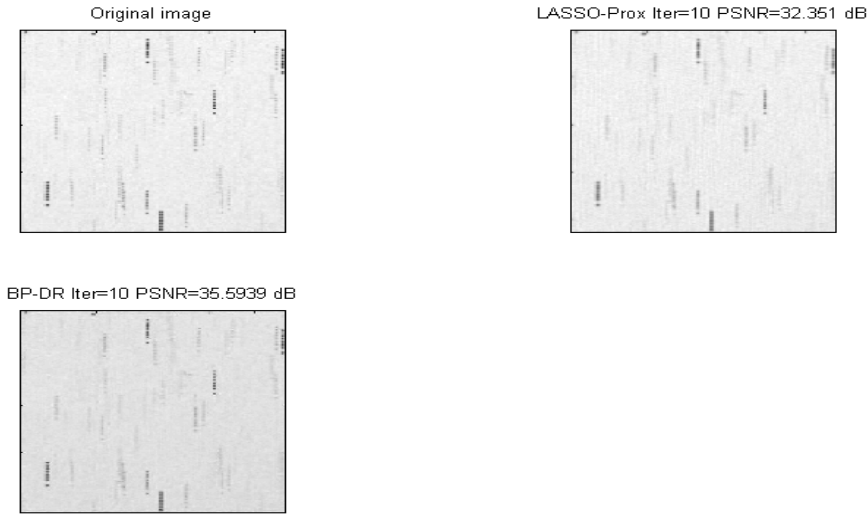

Figure 7

We see again that there is no big difference between the two methods.

### 4.2.2 Image from the Rozhen plate archive

We have chosen an image of the **Pleiades** from the Rozhen plate archive. The details of the algorithms are as above. We take a **PHSD wavelets** with measurement matrix which is 50 lines in the Fourier domain:

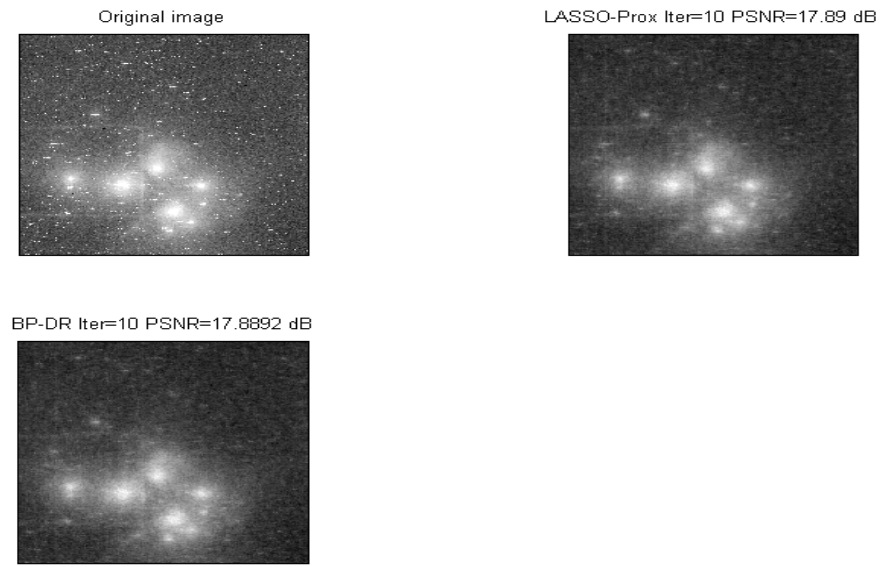

Figure 8

We compare the same with the **Daubechies wavelets**:

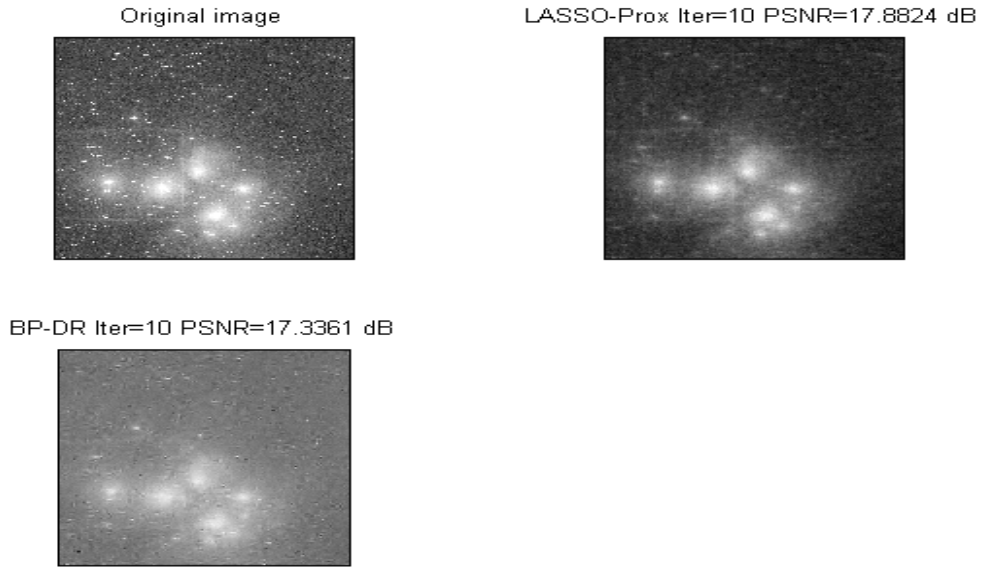

Figure 9

We see that there is no essential difference between the two approaches.

If we increase the size of the measurement matrix ( 100 lines in the Fourier domain) we obtain with **PHSD wavelets** the following:

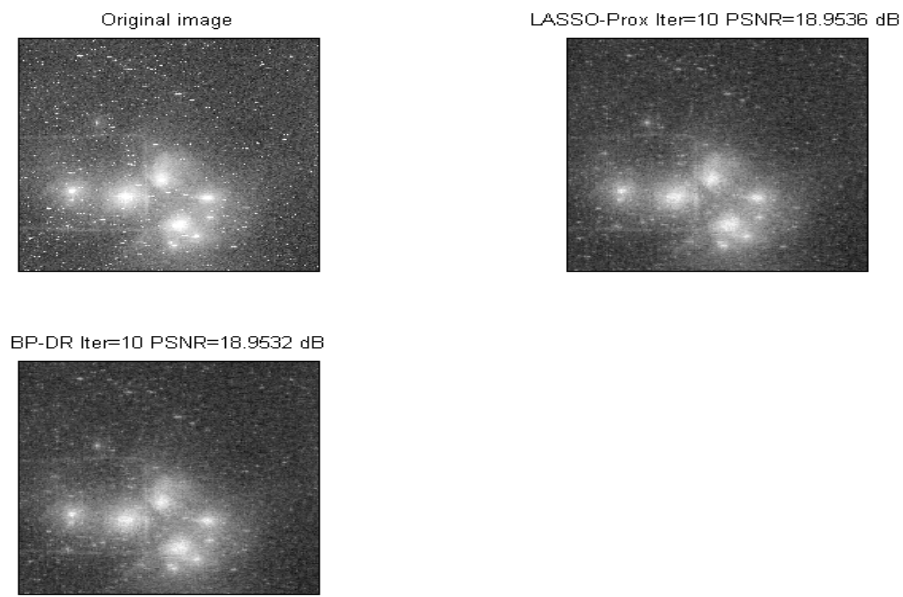

Figure 10

With **Daubechies wavelets** we obtain the following Figure (100 lines in the Fourier domain):

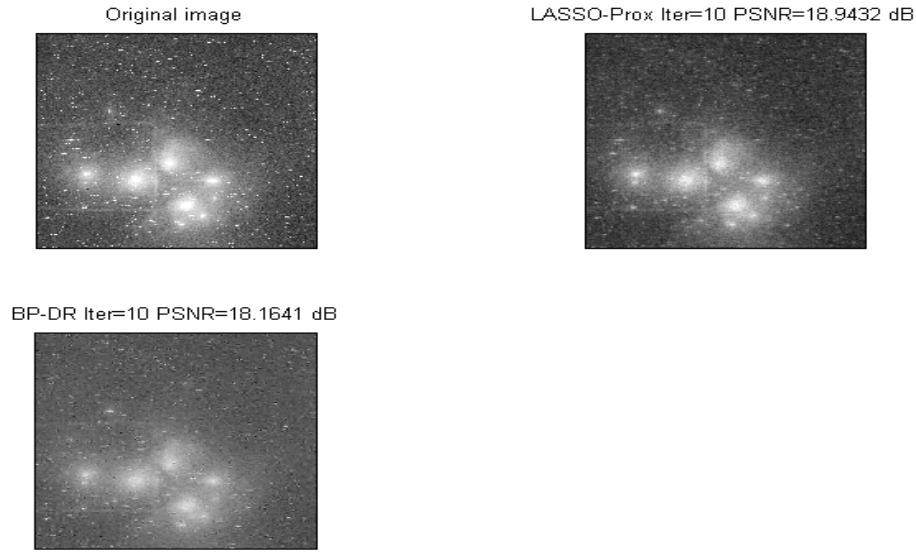

Figure 11

## 5 Final Remarks

In the recent monograph [14] the scheme of CS is applied by using wavelet transforms (Haar wavelets). In [13] the authors apply the CS approach by using the curvelets basis of Candes and Donoho, cf. [5]. In the papers [12], [15] and [8], one finds the Quantization interplay with CS frameworks.

**Acknowledgement.** The first named author was sponsored partially by the Alexander von Humboldt Foundation, and both authors were sponsored by Project DO--2-275/2008 " Astroinformatics" with Bulgarian NSF.

**ABOUT THE AUTHORS:**

Ognyan Kounchev, Prof., Dr., Institute of Mathematics and Informatics, Bulgarian Academy of Science & IZKS, University of Bonn; kounchev@gmx.de

Damyan Kalaglarsky, Institute of Astronomy, Bulgarian Academy of Science; damyan@skyarchive.org.